\documentclass[titlepage,11pt]{article}
\usepackage{latexsym}
\usepackage{amsfonts}
\usepackage{amsmath}
\usepackage[dvips]{graphicx}
\usepackage{color}
\usepackage{epsfig}
\usepackage{epstopdf}
\newcommand\blackslug{\hbox{\hskip 1pt \vrule width 4pt height 8pt depth 1.5pt
        \hskip 1pt}}
\newcommand\bbox{\hfill \quad \blackslug \bigbreak}

\def\l{,\ldots,}

\title{Disjoint dijoins}
\author{Maria Chudnovsky\thanks{Supported by NSF grants IIS-1117631 and DMS-1265803.}\\
Columbia University, New York, NY 10027, USA
\\
\\
Katherine Edwards\thanks{Supported by an NSERC graduate fellowship.}\\
Princeton University, Princeton, NJ 08540, USA
\\
\\
Ringi Kim\\
Princeton University, Princeton, NJ 08544, USA
\\
\\
Alex Scott\\
University of Oxford, Oxford OX1 3LB, UK
\\
\\
Paul Seymour\thanks{Supported by ONR grant N00014-10-1-0680 and NSF grant DMS-1265563.}\\
Princeton University, Princeton, NJ 08544, USA
}

\date{August 23, 2013; revised \today}

\newtheorem{thm}{}[section]

\newcommand{\Proof}{\noindent{\bf Proof.}\ \ }

\begin{document}
\maketitle
\begin{abstract}
A ``dijoin'' in a digraph is a set of edges meeting every directed cut.
D. R. Woodall conjectured in 1976 that if $G$ is a digraph, and every directed cut of
$G$ has at least $k$ edges, then there are $k$ pairwise disjoint dijoins.
This remains open, but 
a capacitated version is known to be false. In particular, 
A. Schrijver gave a digraph $G$ and a subset $S$ of its edge-set, such that 
every directed cut contains at least two edges in $S$, and yet there do not exist
two disjoint dijoins included in $S$. In Schrijver's
example, $G$ is planar, and the subdigraph formed by the edges
in $S$ consists of three disjoint paths.

We conjecture that when $k = 2$,  the disconnectedness of $S$ is crucial: 
more precisely, that if $G$ is a digraph, and $S\subseteq E(G)$ forms
a connected subdigraph (as an undirected graph), and every 
directed cut of $G$ contains at least two edges in $S$, then we can partition
$S$ into two dijoins.

We prove this in two special cases: when $G$ is planar, and
when the subdigraph formed by the edges in $S$ is a subdivision of a 
caterpillar.

\end{abstract}

\section{Introduction}
Some points of terminology, before we begin: in this paper, a {\em graph} 
$G$ consists of a finite set $V(G)$ of {\em vertices}, a finite set $E(G)$ of {\em edges}, 
and an incidence relation
between them; each edge is incident with one or two vertices (its {\em ends}.) A {\em directing} 
of a graph $G$ is a function $\eta$ with domain $E(G)$, where $\eta(e)$ is an end of $e$ for each $e\in E(G)$ (we call
$\eta(e)$ the {\em head} of $e$.) A {\em digraph} $G$ consists of a graph (denoted by $G^-$) and a directing of $G^-$.
If $e$ is an edge of a graph with ends $u,v$, we sometimes refer to the ``edge $uv$''. If $G$ is a digraph, and we refer to an
edge $uv$, this means ``the edge $uv$ of $G^-$'', and does {\em not} imply that this edge has head $v$.
Our other definitions are standard.

Let $G$ be a digraph.
If $X\subseteq V(G)$, $D^+(X)=D_G^+(X)$ denotes the set of all edges of $G$ with tail in $X$ and head in $V(G)\setminus X$,
and $D^-(X) = D^+(V(G)\setminus X)$.
A {\em directed cut} of $G$ means a set $C$ of edges such that there exists $X\subseteq V(G)$ with 
$X, V(G)\setminus X\ne \emptyset$, and $D^-(X) = \emptyset$ and $D^+(X) = C$.
A {\em dijoin} means a subset of $E(G)$ with nonempty intersection with every directed cut of $G$.
D.R.Woodall~\cite{woodall} proposed the following conjecture in 1976:

\begin{thm}\label{woodallconj}
{\bf Woodall's conjecture.}
Let $G$ be a digraph and $k\ge 0$ an integer such that every directed cut has at least $k$ edges. Then there are $k$
pairwise disjoint dijoins.
\end{thm}

This is easily proved for $k\le 2$, but it is still open for $k = 3$, even for planar digraphs $G$.
Interest in \ref{woodallconj} stems from the Lucchesi-Younger theorem~\cite{lucchesi}, which is in some sense dual:

\begin{thm}\label{lucchesiyounger}
Let $G$ be a digraph and $k\ge 0$ an integer such that every dijoin has at least $k$ edges. Then there are $k$
pairwise disjoint directed cuts.
\end{thm}

The Lucchesi-Younger theorem remains true in a capacitated version, as follows ($\mathbb{Z}^+$ denotes the set of non-negative integers):
\begin{thm}\label{LYcap}
Let $G$ be a digraph, and $c$ a map from $E(G)$ to $\mathbb{Z}^+$,  and $k\ge 0$ an integer such that $\sum_{e\in D}c(e)\ge k$
for every dijoin $D$. Then there are $k$
directed cuts such that every edge $e$ is in at most $c(e)$ of them.
\end{thm}
This is easily deduced from \ref{lucchesiyounger} by replacing every edge $e$ with $c(e)>0$ by a directed path of length $c(e)$, 
and by contracting every edge $e$ with $c(e)=0$.

However, the corresponding capacitated version of \ref{woodallconj} is false; indeed, it is false even if $c$ is $0,1$-valued
and $k=2$ and $G$ is planar, as an example due to A.Schrijver~\cite{schrijver} shows. (Figure 1.)
In this paper we investigate further the case when $c$ is $0,1$-valued
and $k=2$.

If $S,T$ are graphs or digraphs, we say they are {\em compatible} if they have the same vertex set
and $E(S)\cap E(T) = \emptyset$. Thus if $S,T$ are both graphs, then they are both subgraphs of a graph $S\cup T$, and if they
are both digraphs then similarly $S\cup T$ is a digraph.
Our problem in this paper is: let 
$S,T$ be compatible digraphs, such that every directed cut of $S\cup T$ has at least two edges in $S$.
When does it follow that $S$ can be partitioned into two dijoins of $S\cup T$?
Schrijver's example shows that this is not always true, even if $S\cup T$ is planar. 
\begin{figure} [h!]
\centering
\includegraphics{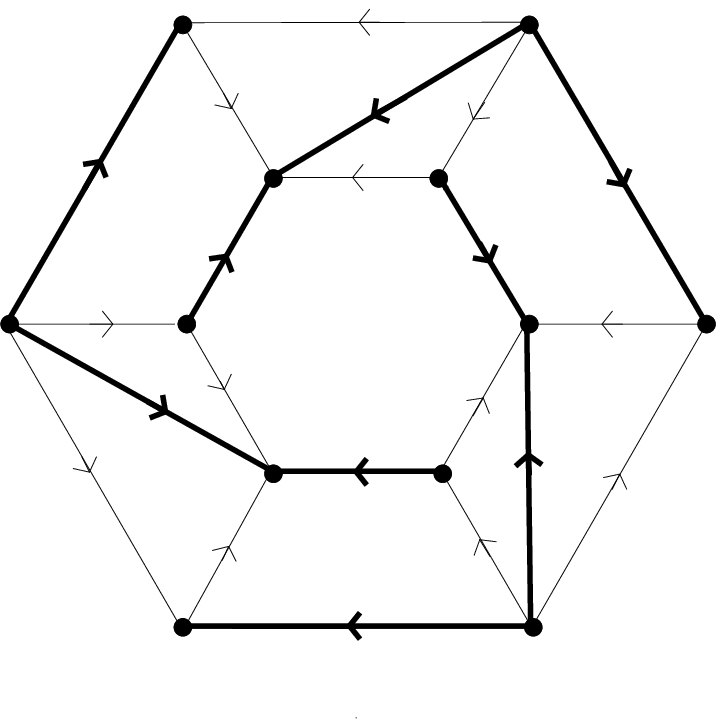}
\caption{Schrijver's counterexample. $S$ contains the thick edges.}
\label{Schrijver}
\end{figure}

When $S,T$ are as in Schrijver's counterexample, 
$S^-$ consists of three vertex-disjoint paths, each of length three. Two more counterexamples were
given in a paper
of Cornu\'ejols and Guenin~\cite{cornuejols}, and Williams~\cite{williams} gave more, all of them minimal 
(in some sense that we do not discuss here).

But all these counterexamples have the property that $S^-$ consists of three disjoint paths; and in this paper we ask, what if $S^-$
is connected? We will study the following conjecture:

\begin{thm}\label{k=2conj}
{\bf Conjecture:} Let $S,T$ be compatible digraphs such that $S^-$ is connected and every directed cut of $S\cup T$ has 
at least two edges in $S$. Then $E(S)$ includes two disjoint dijoins of $S\cup T$.
\end{thm}

Perhaps one can extend this as follows:
\begin{thm}\label{kgenconj}
{\bf Conjecture:} Let $S,T$ be compatible digraphs such that $S^-$ is connected and every directed cut of $S\cup T$ has
at least $k$ edges in $S$. Then $E(S)$ includes $k$ pairwise disjoint dijoins of $S\cup T$.
\end{thm}
This evidently implies Woodall's conjecture, and although it seems much too strong, we have failed to
disprove it so far. (As far as we know, it might be true even if we allow $S^-$ to have two components.)
But in this paper we have nothing more to say about \ref{kgenconj}, and will confine ourselves to 
\ref{k=2conj}.

We have two main results, proofs of two special cases of \ref{k=2conj}, the following. 
Let us say that a tree $T$ is a {\em caterpillar subdivision} if there is a path $P$
of $T$ containing every vertex of $T$ with degree at least three.
\begin{thm}\label{mainthm}
Let $S,T$ be compatible digraphs, such that $S^-$ is connected and every directed cut of $S\cup T$ has
at least two edges in $S$. Suppose that either
\begin{itemize}
\item $S^-$ is a caterpillar subdivision, or
\item $S\cup T$ is planar.
\end{itemize}
Then $E(S)$ can be partitioned into two dijoins.
\end{thm}

The second result is particularly pleasing
because Schrijver's counterexample is planar.
We prove the first assertion of \ref{mainthm} in section 3, and the second in section 4.

\section{Orienting a tree}

It is convenient to work instead with a modified form of \ref{k=2conj}, the following. If $G$ is a digraph, we
say that $X\subseteq V(G)$ is an {\em outset} of $G$ if $D^-_G(X)=\emptyset$ and $X\ne \emptyset, V(G)$.
If $G$ is a graph or digraph and $X\subseteq V(G)$, $D(X)$ or $D_G(X)$ denotes the set of all edges of $G$ with
an end in $X$ and an end in $V(G)\setminus X$.

\begin{thm}\label{mainconj}
{\bf Conjecture: }Let $S,T$ be respectively a tree and a digraph, compatible. Suppose that $|D_S(X)|\ge 2$
for every outset $X$ of $T$. Then there is a directing of $S$, forming a digraph $S'$ say, such that 
$D^+_{S'}(X),D^-_{S'}(X)\ne \emptyset$ for every outset $X$ of $T$.
\end{thm}

(Yes, the tree is $S$ and not $T$, but soon we will take planar duals and $T$ will become the tree, which
is why we chose this notation.)
The statement of \ref{mainconj} makes sense if $S$ is permitted to be something different from a tree, and it is easy to
see that cycles in $S$ can be handled just by contracting them, so we could assume that $S$ is a forest. 
But it is false when $S$ is a forest with three components; Schrijver's counterexample (figure 1) is still a counterexample 
(remove the directions from the edges in $S$.)

It is easy to see (we do not need the result, so we omit the proof) that \ref{mainconj} is true for all pairs
$S,T$ if and only if \ref{k=2conj} is true for all pairs $S,T$. (The proof that \ref{mainconj} implies
\ref{k=2conj} is like the proof below
that \ref{orientthm} implies \ref{mainthm}. For the converse, consider replacing each edge of $S$ by a path of two
oppositely directed edges.) We have checked (on a computer) that \ref{mainconj} 
is true for all trees $S$ with at most twelve vertices of degree different from two (and any number of degree two), 
but in this paper we prove the following two statements:

\begin{thm}\label{orientthm} 
Let $S,T$ be respectively a tree and a digraph, compatible. Suppose that $|D_S(X)|\ge 2$
for every outset $X$ of $T$. Suppose in addition that either
\begin{itemize}
\item $S$ is a caterpillar subdivision, or
\item $S\cup T^-$ is planar.
\end{itemize}
Then there is a directing of $S$, forming a digraph $S'$ say, such that
$D^+_{S'}(X),D^-_{S'}(X)\ne \emptyset$ for every outset $X$ of $G$.
\end{thm}

Before we prove \ref{orientthm}, let us show that it implies \ref{mainthm}.
If $G$ is a graph
or digraph, and $X\subseteq E(G)$,
$G/X$ denotes the graph or digraph obtained by contracting all the edges in $X$, and when $e$ is an edge we write $G/e$ for $G/\{e\}$.

\bigskip

\noindent{\bf Proof of \ref{mainthm}, assuming \ref{orientthm}.}
Let $(S,T)$ be as in \ref{mainthm}, and let us prove \ref{mainthm} by induction on $|V(S)|$. Let $G = S\cup T$.
Suppose first that the graph
$S^-$ has a cycle $C$ (and so $S^-$ is not a tree, and therefore $S\cup T$ is planar). Let $G_1$ be obtained from the digraph
$G$ by contracting all edges in $E(C)$, and let $S_1,T_1$ be the subdigraphs of $G_1$ with vertex sets $V(G_1)$ and with edge
sets $E(G_1)\cap E(S), E(G_1)\cap E(T)$ respectively.
Then $S_1,T_1$ are compatible, and $G_1 =S_1\cup T_1$ is planar, and every directed cut of $G_1$ is a directed cut of $G$ and hence
has at least two edges in $S_1$. From the inductive hypothesis, $E(S_1)$ can be partitioned into two dijoins $A_1,B_1$ of $G_1$.
Choose an orientation of $C$, and let $A_2$ be the set of edges of $C$ that are positively oriented, and $B_2 = E(C)\setminus A_2$.
Every directed cut of $S\cup T$ that is not a directed cut of $G_1$ contains an edge in $E(C)$, and hence contains an edge in $A_2$
and an edge in $B_2$. It follows that $A_1\cup A_2, B_1\cup B_2$ are dijoins of $G$, as required.

Thus we may assume that $S^-$ is a tree. Suppose that
there is an outset $X$ of $T$ such that $|D_S(X)|\le 1$. Every directed cut of $G$ contains at least two edges in $S$;
so $D_G(X)$ is not a directed cut of $G$. Since $X$ is an outset of $T$, it follows that $|D_S(X)| = 1$, say 
$D_S(X) \cap S = \{s\}$, 
and $s\in D_S^-(X)$.
Let $G_1 = G/s$ and let $S_1,T_1$ be the subdigraphs of $G_1$ with vertex sets $V(G_1)$ and  with edge
sets $E(G_1)\cap E(S), E(G_1)\cap E(T)$ respectively.
Thus if $G$ is planar then so is $G_1$, and if $S^-$
is a caterpillar subdivision then so is $S_1^-$. From the inductive hypothesis, $E(S_1)$ can be partitioned into two dijoins
$A,B$ of $G_1$. We claim that $A,B$ are dijoins of $G$. For suppose that $A$ is not, say; then there is an outset $Y$ of $G$
such that $A\cap D_G(Y) = \emptyset$. Since $A$ is a dijoin of $G_1$, it follows that $s\in D_G(Y)$, and hence
$s\in D^+_G(Y)$ since $Y$ is an outset. 
We recall that $D_S^-(X) = \{s\}$, and $D_S^+(X)=\emptyset$, and since $A\subseteq E(S)$ it follows that
$A\cap D(X) = \emptyset$. Let $s$ have tail $u$ and head $v$.
Thus $u\in Y\setminus X$ and $v\in X\setminus Y$. Now $D_G^-(X\cap Y) = \emptyset$; because there are no edges of $G$ 
from $Y\setminus X$
to $X\cap Y$ since $D_G^-(X) = \{s\}$, and none from $V(G)\setminus Y$ to $Y$ since $D_G^-(Y) = \emptyset$.
Since $A$ is a dijoin of $G_1$ and $A\cap D^-(X\cap Y) =\emptyset$, it follows that $D^-(X\cap Y)$ is not a directed cut
of $G_1$; and so $X\cap Y = \emptyset$. Similarly $X\cup Y = V(G)$, and so $X=V(G)\setminus Y$. But $D(X)$ contains only one 
edge in $S$, and yet $D(Y)$ is a directed cut of $G$ and so contains at least two edges of $S$ by hypothesis, a contradiction.
This proves that $A,B$ are dijoins of $G$, and hence the result holds.

We may therefore assume that $|D_S(X)|\ge 2$
for every outset $X$ of $T$. Since either $S^-\cup T^-$ is planar or $S^-$ is a caterpillar subdivision,
\ref{orientthm} implies that there is a directing of $S^-$, forming a digraph $S'$ say, such that
$D^+_{S'}(X),D^-_{S'}(X)\ne \emptyset$ for every outset $X$ of $T$. Let $A$ be the set of edges in $S$ that have 
the same head in $S$ and in $S'$, and $B$ those given different heads. We claim that $A,B$ are dijoins of $G$.
For let $X$ be an outset of $G$. Thus $X$ is an outset of $T$, and
so $D^+_{S'}(X),D^-_{S'}(X)\ne \emptyset$. But $D^+_{S'}(X)\subseteq A$, and $D^-_{S'}(X)\subseteq B$, and so
both $A,B$ have nonempty intersection with $D_G(X)$. This proves that $A,B$ are dijoins, and completes the proof of \ref{mainthm}.~\bbox

\section{Caterpillars and biases}

In this section we prove the first assertion of \ref{orientthm}.
Let $G$ be a graph or digraph. If $A,B\subseteq V(G)$, we denote by $D(A,B)$ or $D_G(A,B)$ the set of all edges of $G$ with an end in $A\setminus B$
and an end in $B\setminus A$. If $S$ is a graph, a {\em bias} in $S$ is a set $\mathcal{B}$ of subsets of $V(S)$ such that
\begin{itemize}
\item if $A\in \mathcal{B}$ then $A\ne \emptyset, V(S)$
\item if $A,B\in \mathcal{B}$ and $D(A,B) = \emptyset$ then $A\cap B, A\cup B\in \mathcal{B}\cup \{\emptyset, V(S)\}$;
\item if $A,B\in \mathcal{B}$ and $|D(A,B)| = 1$ then at least one of $A\cap B, A\cup B$ belongs to $\mathcal{B}\cup \{\emptyset, V(S)\}$.
\end{itemize}
We propose the following conjecture.
\begin{thm}\label{biasconj}
Let $S$ be a tree, and let $\mathcal{B}$ be a bias in $S$, such that $|D_S(X)|\ge 2$ for each $X\in \mathcal{B}$.
Then there is a directing of $S$, forming a digraph $S'$ say, such that
$D^+_{S'}(X),D^-_{S'}(X)\ne\emptyset$ for each $X\in \mathcal{B}$.
\end{thm}

This conjecture implies \ref{mainconj}; because if $S,T$ are as in \ref{mainconj} then
the set of all outsets of $T$ 
is a bias $\mathcal{B}$ say in $S$, 
and then the truth of \ref{k=2conj} for $S,T$ is implied
by the truth of \ref{biasconj} for $S,\mathcal{B}$. We have not been able to prove \ref{biasconj} in general, but we have
checked it on a computer for all trees $S$ with at most twelve vertices.

Let us say a forest $S$ is {\em upright} if \ref{biasconj} holds for every bias; that is,
for every bias $\mathcal{B}$ in $S$ such that $|D_S(X)|\ge 2$ for each $X\in \mathcal{B}$,
there is a directing of $S$, forming a digraph $S'$ say, such that
$D^+_{S'}(X),D^-_{S'}(X)\ne\emptyset$ for each $X\in \mathcal{B}$.
Thus, conjecture \ref{biasconj} says that every tree is upright. We do not know whether every forest with two components is upright;
but three components (all with edges) is too many. Even the forest with six vertices and three pairwise disjoint edges is not upright.
To see this, let $S$ have edges $v_1v_2, v_3v_4,v_5v_6$ say; then 
$$\{\{v_1,v_3,v_5\},\{v_1,v_4,v_6\},\{v_2,v_3,v_6\},\{v_2,v_4,v_5\}\}$$
is a bias, and $S$ cannot be oriented to satisfy the conclusion of \ref{biasconj}. Perhaps every pair of a forest and a bias
satisfying the hypothesis and not the conclusion of \ref{biasconj} can be contracted to this six-vertex example, which would be a pleasing explanation
why all the counterexamples previously mentioned (due to Schrijver, Cornu\'ejols, Guenin, Williams) have at least 
three components all with edges.

Our task in this section is to show that every caterpillar subdivision is upright, thereby proving the first statement of \ref{orientthm}.
Note that if $\mathcal{B}$ is a bias in $S$, then so is the set of all sets $V(S)\setminus X\;(X\in \mathcal{B})$, and we call the
latter the {\em reverse bias} of $\mathcal{B}$.

The advantage of biases is the following convenient lemma:
\begin{thm}\label{series}
Let $S$ be a tree, and let $v$ be a vertex of $S$ with degree two. Let $u,w$ be its neighbours, and let
$T$ be the tree obtained from $S$ by deleting $v$ and adding a new edge joining $u,w$. If $T$ is upright then so is $S$.
\end{thm}
\Proof
Let $\mathcal{B}$ be a bias in $S$, such that $|D_S(X)|\ge 2$ for each $X\in \mathcal{B}$.
We say that a subset $X\subseteq V(S)$ is {\em linear} if its intersection with $\{u,v,w\}$ is one of
$$\emptyset, \{u\},\{u,v\},\{u,v,w\},\{v,w\},\{w\}.$$
Let $\mathcal{C}$ be the set of all $Y\subseteq V(T)$ such that there is a linear $X\in \mathcal{B}$ with
$X\cap V(T) = Y$. (We call such a set $X$ a {\em parent} for $Y$.)
\\
\\
(1) {\em $\mathcal{C}$ is a bias in $T$.}
\\
\\
For let $Y\in \mathcal{C}$, and let $X$ be a parent for $Y$.
If $Y=\emptyset$, then $u,w\notin X$, and so $v\notin X$ since $X$ is linear; and so $X=\emptyset$,
a contradiction. Thus $Y\ne \emptyset$, and similarly $Y\ne V(T)$.

Now let $Y_1,Y_2\in \mathcal{C}$, and for $i = 1,2$ let $X_i$ be a parent of $Y_i$. Let $|D_T(Y_1,Y_2)| = n$,
where $n\le 1$. Let $m = |D_S(X_1,X_2)|$; 
it follows that $m\le n$, since $X_1,X_2$ are linear. On the other hand, since $X_1,X_2$ are linear, it follows
that at least one of $X_1\cap X_2, X_1\cup X_2$ is linear. But
at least $2-m$ of $X_1\cup X_2,X_1\cap X_2$ belong to $\mathcal{B}\cup \{\emptyset, V(S)\}$.
If at least $2-n$ of these sets are linear, then it follows that 
$2-n$ of $Y_1\cup Y_2,Y_1\cap Y_2$ belong to $\mathcal{C}\cup \{\emptyset, V(T)\}$ as required.
We may therefore assume that at least $n-m+1$ of $X_1\cup X_2,X_1\cap X_2$ belong to $\mathcal{B}\cup \{\emptyset, V(S)\}$ 
and are not linear. Since as we saw, at most one of them is not linear, we deduce that $n-m+1\le 1$, and since $m\le n$ it follows that $m=n$.
But since $X_1,X_2$ are linear and one of $X_1\cap X_2, X_1\cup X_2$ is not, it follows that the edge $uw$ of $T$ belongs to 
$D_T(Y_1,Y_2)$, and neither of the edges $uv,vw$ of $S$ belong to $D_S(X_1,X_2)$, contradicting that $m = n$.
This proves (1).
\\
\\
(2) {\em $|D_T(Y)|\ge 2$ for each $Y\in \mathcal{C}$.}
\\
\\
For let $X$ be a parent of $Y$. Then $|D_S(X)|\ge 2$ by hypothesis. We claim that $|D_T(Y)|\ge |D_S(X)|$. If $D_S(X)\subseteq D_T(Y)$
this is true, and so
we may assume that one of the edges
$uv,vw$ belongs to $D_S(X)$. Since $X$ is linear, it follows that 
only one of $uv,vw$ belongs to $D_S(X)$, and also $uw$
belongs to $D_T(Y)$; and consequently $|D_T(Y)|\ge |D_S(X)|$. This proves (2).

\bigskip

Since $T$ is upright, there is a directing of $T$, forming a digraph $T'$, such that 
$D^+_{T'}(Y),D^-_{T'}(Y)\ne\emptyset$ for each $Y\in \mathcal{C}$. Now we define a directing $\eta$ of $S$ as follows.
For each edge $e$ of $S$ different from $uv,vw$
let $\eta(e)$ be the head of $e$ in $T'$. For $uv,vw$, choose $\eta(uv),\eta(vw)$ so that one of them is $v$ and one of them is the
head of $uw$ in $T'$. 
We claim that $D^+_{S'}(X),D^-_{S'}(X)\ne\emptyset$ for each $X\in \mathcal{B}$.
For if $X$ is not linear, then one of $uv,vw$ is in $D^+_{S'}(X)$ and the other is in $D^-_{S'}(X)$, so both these sets are nonempty.
If $X$ is linear, we can assume by taking the reverse bias if necessary that $v\notin X$, and so $X\in \mathcal{C}$; but then
$D^+_{T'}(X),D^-_{T'}(X)\ne\emptyset$, and so $D^+_{S'}(X),D^-_{S'}(X)\ne\emptyset$ as required. This proves \ref{series}.~\bbox

A path $P$ of a tree $S$ is a {\em spine} if every vertex of $S$ either belongs to $P$ or has a  
neighbour in $P$; and a tree that has a spine is called a {\em caterpillar}.
In view of \ref{series}, in order to prove the first assertion of \ref{orientthm} it suffices to prove the following.

\begin{thm}\label{caterpillar}
Let $S$ be a caterpillar; then $S$ is upright. Moreover, let $P$ be a spine of $S$,
and let $\mathcal{B}$ be a bias in $S$, such that
$|D_S(X)|\ge 2$ for each $X\in \mathcal{B}$.
Then there is a directing of $S$, forming a digraph $S'$ say, such that
$D^+_{S'}(X),D^-_{S'}(X)\ne\emptyset$ for each $X\in \mathcal{B}$, and such that $P$ corresponds to a directed path of $S'$.
\end{thm}
\Proof
The second assertion implies the first, and we prove the second by induction on $|V(S)|$. Thus, let
$P$ have vertices $p_1\l p_n$ in order, where $n\ge 1$. For $1\le i\le n$, let $Q_i$ be the set of vertices of $S$ not in $V(P)$
that are adjacent to $p_i$. By choosing $P$ maximal, we may assume that $Q_1,Q_n$ are empty. If $n=1$ the claim is trivial,
so we assume that $n\ge 2$. Let $T= S\setminus\{p_1\}$. For $Y\subseteq V(T)$, let $\tilde{Y} =Y$ if $p_2\notin Y$, 
and $\tilde{Y} = Y\cup \{p_1\}$ if $p_2\in Y$. Let $\mathcal{C}$ be the set of all subsets $Y\subseteq V(T)$
such that $\tilde{Y}\in \mathcal{B}$.
\\
\\
(1) {\em $\mathcal{C}$ is a bias in $T$, and $|D_T(Y)|\ge 2$ for each $Y\in \mathcal{C}$.}
\\
\\
For let $Y\in \mathcal{C}$. Since $\tilde{Y}\in \mathcal{B}$, it follows that $\tilde{Y}\ne \emptyset, V(S)$, and so 
$Y\ne \emptyset, V(T)$. Now let $Y_1,Y_2\in \mathcal{C}$, with $|D_T(Y_1,Y_2)|\le 1$. Then $|D_S(\tilde{Y_1},\tilde{Y_2})|\le 1$,
and so one of $\tilde{Y_1}\cup \tilde{Y_2},\tilde{Y_1}\cap \tilde{Y_2}\in \mathcal{B}\cup\{\emptyset, V(S)\}$, and both if 
$D_S(\tilde{Y_1},\tilde{Y_2})=\emptyset$. Hence one of $Y_1\cup Y_2,Y_1\cap Y_2$ is in $\mathcal{C}\cup\{\emptyset, V(T)\}$, and both
if $D_T(Y_1,Y_2)=\emptyset$. This proves that $\mathcal{C}$ is a bias in $T$. But if $Y\in \mathcal{C}$, then $D_T(Y) = D_S(\tilde{Y}))$, and so
$|D_T(Y)|\ge 2$. This proves (1).
\\
\\
(2) {\em There is a directing of $S$, forming a digraph $S'$ say, such that 
\begin{itemize}
\item $p_ip_{i+1}$ has head $p_{i+1}$ for $1\le i<n$, and so $P$ becomes a directed path of $S'$, directed from $p_1$ to $p_n$, and
\item $D^+_{S'}(\tilde{Y}),D^-_{S'}(\tilde{Y})\ne\emptyset$ for each $Y\in \mathcal{C}$.
\end{itemize}
}
\noindent Since the path with vertices $p_2\l p_n$ in order is a spine of $T$, (1) implies that there is a directing of $T$, forming
a digraph $T'$, such that $D^+_{T'}(Y),D^-_{T'}(Y)\ne\emptyset$ for each $Y\in \mathcal{C}$, and such that $P\setminus p_1$ 
becomes a directed path of $T'$. By reversing the direction of all edges if necessary, we may assume that
$p_2$ is the first vertex of this directed path; and now extend the directing to one of $S$ by assigning $p_2$ to be the 
head of $p_1p_2$. Now let $Y\in \mathcal{C}$. Then $D^+_{T'}(Y),D^-_{T'}(Y)\ne\emptyset$, but $D^+_{S'}(\tilde{Y})=D^+_{T'}(Y)$
and $D^-_{S'}(\tilde{Y})=D^-_{T'}(Y)$, and so $D^+_{S'}(\tilde{Y}),D^-_{S'}(\tilde{Y})\ne\emptyset$. 
This proves (2).

\bigskip
Among all choices of $S'$ satisfying (2), choose one such that as many edges of $S$ as possible are directed such that
their head is not in $\{p_1\l p_n\}$.
\\
\\
(3) {\em There is no $X\in \mathcal{B}$ such that $p_1\in X$, $p_2\notin X$, and
$D^-_{S'}(X)= \emptyset$.}
\\
\\
For suppose that such a set $X$ exists. Choose $X$ such that $|X|$ is minimum.
Since $X$ is an outset of $S'$, and since $p_2\notin X$, it follows that $p_3\l p_n\notin X$. Since $X\in \mathcal{B}$ and therefore
$|D_S(X)|\ge 2$, it follows that there exists $i\in \{2\l n-1\}$ such that $X\cap Q_i\ne \emptyset$. Choose such a value of
$i$, maximum, and let $v\in X\cap Q_i$. (See figure 2.)
\begin{figure} [h!]
\centering
\includegraphics{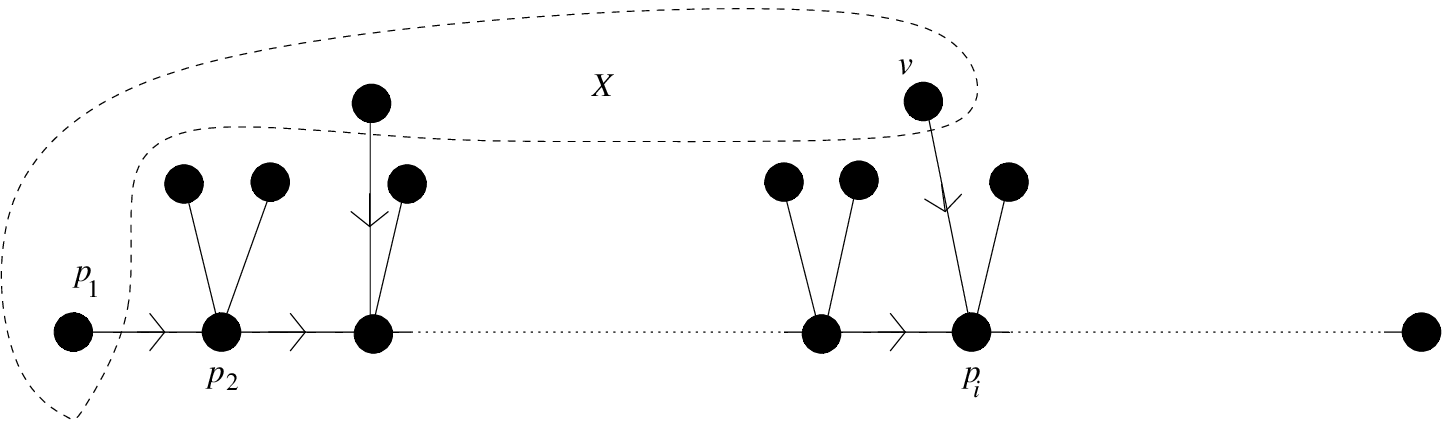}
\caption{}
\label{catfig}
\end{figure}
Since $X$ is an outset of $S'$, the edge $e=p_iv$ of $S$ has head $p_i$ in $S'$.
From the choice of $S'$, reversing the direction of this edge does not give a better choice of $S'$ satisfying (2), and so there
exists $Y\in \mathcal{C}$ such that one of $D^+_{S'}(\tilde{Y}),D^-_{S'}(\tilde{Y})= \{e\}$. Thus either
$v\in Y,p_i\notin Y$, and $D^+_{S'}(\tilde{Y})=\{e\}$, or $p_i\in Y, v\notin Y$ and $D^-_{S'}(\tilde{Y})=\{e\}$.

Suppose first that $v\in Y,p_i\notin Y$, and $D^+_{S'}(\tilde{Y})=\{e\}$. Since $p_i\notin Y$ it follows that $p_1\l p_{i-1}\notin Y$.
Consequently $D_S(X,\tilde{Y}) = \emptyset$, and since $\mathcal{B}$ is a bias, it follows that 
$X\cap \tilde{Y}\in \mathcal{B}$. By hypothesis, $|D_S(X\cap \tilde{Y})|\ge 2$, and so there exists $u\in X\cap \tilde{Y}$
different from $v$. Let $u\in Q_h$ say; then $h\le i$ from the choice of $i$, and so $p_h\notin Y$. Since
$X$ is an outset of $S'$, it follows that $u$ is the tail of the edge $up_h$, contradicting that $D^+_{S'}(\tilde{Y})=\{e\}$.
Thus $p_i\in Y, v\notin Y$ and $D^-_{S'}(\tilde{Y})=\{e\}$. It follows that $p_1\l p_{i-1}\in \tilde{Y}$. 

We claim that $X\subseteq \tilde{Y}\cup \{v\}$. For suppose that there exists $u\ne v$ with $u\in X$ and 
$u\notin \tilde{Y}$. Since $p_1\in \tilde{Y}$ it follows that $u\in Q_h$ for some $h$ with $2\le h\le i$; let $f$ be the edge $up_h$.
Since $X$ is an outset of
$S'$, $p_h$ is the head of $f$ in $S'$, contradicting that $D^-_{S'}(\tilde{Y})=\{e\}$. This proves that 
$X\subseteq \tilde{Y}\cup \{v\}$, and consequently $|D_S(X,\tilde{Y})| = 1$.

Since $\mathcal{B}$ is a bias, one of 
$X\cup \tilde{Y},X\cap \tilde{Y}\in \mathcal{B}\cup\{\emptyset, V(S)\}$. We assume first that 
$X\cup \tilde{Y}\in \mathcal{B}\cup\{\emptyset, V(S)\}$.
Since $X\cup \tilde{Y} = \tilde{Y}\cup \{v\}$ and $D^-_{S'}(\tilde{Y}\cup \{v\})=\emptyset$, it follows that 
$Y\cup \{v\}\in \mathcal{C}\cup \{\emptyset, V(T)\}$ and $D^-_{T'}(Y\cup \{v\})=\emptyset$. From the choice of $T'$, it follows that
$Y\cup \{v\}= V(T)$; but then $|D_{S}(\tilde{Y})|=1$, contrary to hypothesis. Finally, we assume that
$X\cap \tilde{Y}\in \mathcal{B}\cup\{\emptyset, V(S)\}$. But $X\cap \tilde{Y} = X\setminus \{v\}$, and 
since $p_1\in X\cap \tilde{Y}$, it follows that $X\setminus \{v\} \in \mathcal{B}$, and yet $D^-_{S'}(X\setminus \{v\})= \emptyset$,
contrary to the minimality of $X$. This proves (3).
\\
\\
(4) {\em There is no $X\in \mathcal{B}$ such that $p_2\in X$, $p_1\notin X$, and
$D^+_{S'}(X)= \emptyset$.}
\\
\\
This follows by applying (3) to the reverse bias of $\mathcal{B}$ and to the complement of $X$.
\bigskip

Now to complete the proof, we must show that 
if $X\in \mathcal{B}$ then $D^+_{S'}(X),D^-_{S'}(X)\ne  \emptyset$.
This is true from the choice of $T'$ if $X$ contains both or neither of $p_1,p_2$, so we may assume that
$X$ contains exactly one of $p_1,p_2$. If $p_1\in X$ then $D^+_{S'}(X)\ne \emptyset$ since it contains the edge $p_1p_2$,
and $D^-_{S'}(X)\ne \emptyset$ by (3). If $p_2\in X$ and $p_1\notin X$, then the claim follows from (4), similarly.
This proves \ref{caterpillar}, and hence proves the first statement of \ref{orientthm}.~\bbox

\section{Planarity}

In this section we prove the second assertion of \ref{orientthm}, and it is helpful to reword it first.
If $G$ is a digraph, $\overleftarrow{G}$ denotes the digraph obtained by reversing 
the direction of all edges of $G$. In \ref{orientthm}, since $S$ is a spanning tree of $S\cup T^-$,
the statement that there is no outset $X$ of $T$ such that $|D_S(X)|\le 1$ is the same as saying that for each edge $e$ of the
tree $S$, the cutset in $S\cup T^-$ that contains no edge of $S$ except $e$ contains edges of $T$ directed in each direction. 
Also, the conclusion of \ref{orientthm} can be reworded to say that ``we can direct the edges of $S$, forming a digraph $S'$,
such that both $S'\cup T$ and $\overleftarrow{S'}\cup T$ have no directed cuts''.
Finally, if $S,T$ are as in the second assertion of \ref{orientthm},
then since $S$ is a spanning tree of the planar graph $S\cup T^-$, it follows that $E(T)$ is the set of edges of a spanning tree of
the planar dual, and it is convenient to reword everything in terms of the dual. A {\em directed tree} is a digraph $T$ such that
$T^-$ is a tree.  A digraph is {\em acyclic} if it has no directed cycle. We leave the reader to check that after all these modifications, 
the second assertion of \ref{orientthm} becomes the following:

\begin{thm}\label{planarthm}
Let $S$ be a graph and $T$ a directed tree, such that $S,T$ are compatible and $S\cup T^-$ is planar.
Suppose that for every edge $e$ of $S$, the path of $T$ joining the ends of $e$ is not a directed path of $T$. Then
there is a directing of $S$, forming a digraph $S'$ say, such that
the digraphs $S'\cup T, \overleftarrow{S'}\cup T$ are acyclic.
\end{thm}

In contrast with the second assertion of \ref{orientthm},
we can show that the planarity hypothesis in \ref{planarthm} cannot be omitted: let $T$ be the digraph with vertex set $\{v_1\l v_7\}$ and edges the ordered pairs
$$v_1v_2,v_1v_3,v_1v_4,v_5v_2,v_6v_3,v_7v_4,$$
and let $S$ be the graph with the same vertex set and edges the unordered pairs
$$v_1v_5,v_1v_6,v_1v_7,v_2v_6,v_2v_7,v_3v_5,v_3v_7,v_4v_5,v_4v_6.$$
(Figure 3.)
We leave the reader to check that $S$ cannot be directed to satisfy the conclusion of \ref{planarthm}.
\begin{figure} [h!]
\centering
\includegraphics{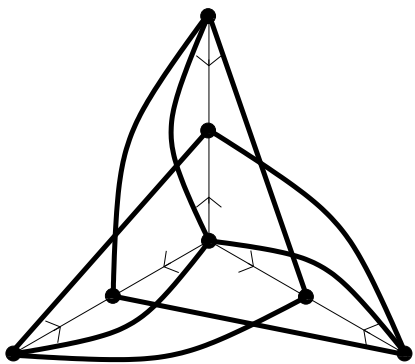}
\caption{Counterexample to the nonplanar extension of \ref{planarthm}. $S$ contains the thick edges.}
\label{nonplanar}
\end{figure}

How else might we try to extend \ref{planarthm}? Let us change the 
hypothesis to say that ``no directed path of $T$ joins the ends of $e$'' so that it makes sense when $T^-$ is not a tree.
It is easy to see that if it is true as stated then it is also true if we just
require that $T^-$ is a forest rather than a tree.
The problems come when $T^-$ is more than a spanning tree rather than less than one; the result is false in general if we permit
$T^-$ to be a theta (the planar dual of Schrijver's counterexample gives a counterexample to this extension of \ref{planarthm}
in which $T^-$ consists of three paths each of length four, with the same ends and otherwise disjoint). We do not know whether
\ref{planarthm} holds if we permit $T^-$ to have exactly one cycle; and in fact we do not know if it holds when $T^-$ is a cycle.

Another way we might try to extend \ref{planarthm} is the following. \ref{planarthm} says there is a directing of $S$ that works for any two
opposite directings of $T$, but perhaps there is a directing of $S$ that works simultaneously for all directings of $T$. 
More exactly, we might hope that:

\bigskip

{\em \noindent Let $S$ be a graph and $T$ a tree, such that $S,T$ are compatible and $S\cup T$ is planar.
Then there is a directing of $S$, forming a digraph $S'$ say, with the following property. Let $T'$ be a digraph with 
$T'^- = T$, such that 
for every edge $e$ of $S$, the path of $T'$ joining the ends of $e$ is not a directed path of $T$. Then
$S'\cup T'$ is acyclic. }

\begin{figure} [h!]
\centering
\includegraphics{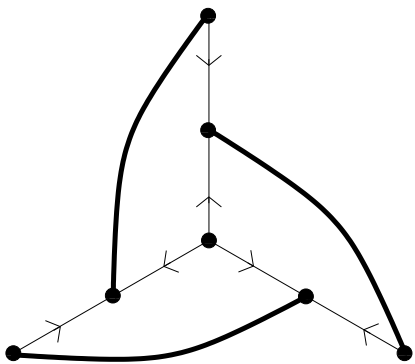}
\caption{Counterexample to an extension of \ref{planarthm}. $S$ contains the thick edges.}
\label{universal}
\end{figure}
\noindent This is true when $T$ is a path, but false when
$T$ is the tree with vertex set $\{v_1\l v_7\}$ and edges the unordered
pairs
$$v_1v_2,v_1v_3,v_1v_4,v_5v_2,v_6v_3,v_7v_4,$$
and $S$ is the graph with the same vertex set and edges the unordered
pairs
$v_2v_6,v_3v_7,v_4v_5$. (Figure 4.)

Let $G$ be a loopless graph drawn in a $2$-sphere, and let $t\in V(T)$. The edges incident with $t$ are drawn in a circular order.
A subset $W\subseteq \{e_1\l e_k\}$ is a {\em $t$-wedge} (with respect to the given drawing)
if it forms an interval of this circular order, that is,
if we can enumerate the edges 
of $G$ incident with $t$ as $f_1\l f_k,f_1$ in circular order, such that $W = \{f_1\l f_i\}$ for some $i$.
 
To prove \ref{planarthm}, it is helpful for inductive purposes to prove a strengthening, the following.

\begin{thm}\label{wedgethm}
Let $S$ be a graph and $T$ a directed tree, such that $S,T$ are compatible, and $S\cup T^-$ is planar.
Suppose that for every edge $e$ of $S$, the path of $T$ joining the ends of $e$ is not a directed path of $T$. 
Fix a drawing of $S\cup T^-$ in a $2$-sphere $\Sigma$,
let $t\in V(T)$, and let $W$ be a $t$-wedge of $S\cup T^-$ with $W\cap E(T) = \emptyset$.
Then
there is a directing of $S$, forming a digraph $S'$ say, such that
\begin{itemize}
\item $S'\cup T$ and $S' \cup \overleftarrow{T}$ are acyclic, and
\item every edge in $W$ has tail $t$.
\end{itemize}
\end{thm}

Why use wedges? There are two other strengthenings that we could try instead,
both of them more natural; we
might try to replace the wedge part of \ref{wedgethm} by
\begin{itemize}
\item let $t$ be a vertex of degree one in $T$; then we can choose $S'$
so that no edge has head $t$, or
\item let $t$ be a vertex of $T$; then we can choose $S'$
so that no edge has head $t$.
\end{itemize}
But neither of these works.
The first is not strong enough; we were unable to make the induction go
through. The second is too strong, because it is false (for a
counterexample, let $T$ have five vertices $v_1\l v_5$, and edges the
ordered pairs
$$v_1v_2,v_3v_2,v_3v_4,v_5v_4,$$
and let $S$ have the same vertex set and edges the unordered pairs
$$v_1v_3, v_1v_4,v_2v_4,v_2v_5,v_3v_5.$$
The wedge form is a compromise between the two that works.

\bigskip

\noindent{Proof of \ref{wedgethm}.\ \ }
We proceed by induction on $|V(T)|$. Let $t$ have degree $d$ in $T$. Thus $d\ge 0$,
and if $d=0$, then $V(T) = \{t\}$, and $E(S) = \emptyset$, and the result is true. So we may assume that $d\ge 1$. 

For simplicity, we use the same notation for a vertex of $S\cup T^-$ and   
its image under the drawing, and the same for an edge and the corresponding open line segment.

For every two vertices
$u,v$ of $T$, $T(u,v)$ denotes the path of $T$ between $u$ and $v$.
Eventually we will choose a directing of $S$, forming a digraph $S'$, and we will need to show that $S'\cup T$ 
and $S'\cup \overleftarrow{T}$ are acyclic. Here is a useful lemma for that purpose. Let $S'$ be a digraph with $S'^- = S$.
A directed cycle $C$ of $S'\cup T$ is {\em optimal} (for $S'$) if $E(C)\cap E(S')$ is minimal. If
$S'\cup T$ is not acyclic then there is a directed cycle in $S'\cup T$ and hence an optimal one.
\\
\\
(1) {\em Let $S'$ be a digraph obtained by directing the edges of $S$, and let $C$ be an optimal directed cycle of  $S'\cup T$.
Then for every directed path $P$ of $T$, $C\cap P$ is either a directed path or null.}
\\
\\
Suppose not; then there are two distinct vertices $u,v$ of $P$, both in $V(C)$ and such that $V(C)$ contains no other 
vertex of the subpath $Q$ of $P$ between $u,v$. Since $P$ is a directed path, so is $Q$. We may assume that
$Q$ has first vertex $u$ and last vertex $v$. Let $P_1,P_2$ be the directed paths in $C$ from $u$ to $v$ and from $v$ to $u$
respectively. Thus $Q\cup P_2$ is a directed cycle of $S'\cup T$, and from the optimality of $C$ it follows that
$E(P_1)\subseteq E(T)$. But then $P_1^-\cup Q^-$ is a cycle of the tree $T^-$, which is impossible. This proves (1).
\\
\\
(2) {\em We may assume that $d\ge 2$.}
\\
\\
Suppose that $d=1$, and let $t_1$ be the neighbour of $t$ in $T$. 
By reversing the directing of $T$ if necessary, we may assume that $t$ is the head of $tt_1$ in $T$.
Let $R$ be the graph obtained as follows. Its vertex set is $V(T)\setminus \{t\}$. For all 
distinct $u,v\in V(R)$, if $u,v$ are adjacent in $S$ then they are adjacent in $R$. In addition, for every edge $e$ of $S$ incident with
$t$ and some other vertex $v$, if $T(v,t_1)$ is not a directed path, then $e$ is an edge of $R$ with ends $v,t_1$.
Note that $R\cup (T\setminus \{t\})^-$ is planar, since it can be obtained from a subgraph of $S\cup T^-$ by contracting
the edge $tt_1$; and for the same reason, there is a $t_1$-wedge $W'$ of the corresponding drawing that contains all 
edges of $E(R)$ that are incident with $t$ in $S$, and contains no edges of $T\setminus t$.
Also, for each edge $e$ of $R$, the path of $T\setminus t$ joining the ends of $e$ is not
a directed path of $T\setminus t$, from the definition of $R$. From the inductive hypothesis, there is an directing of $R$,
forming a digraph $R'$ say, such that $R'\cup (T\setminus t)$ and $R'\cup \overleftarrow{T}\setminus t$ are acyclic, 
and every edge in $E(R)$ that is incident with $t$ in $S$ has head different from $t_1$. 
Define a directing $\eta$ of $S$ as follows, forming a digraph $S'$: 
\begin{itemize}
\item for each edge $e$ of $R$, let $\eta(e)$
be the head of $e$ in $R'$, and 
\item for each edge $e$ of $S$ that is not an edge of $R$, let $\eta(e)$ be its end different from $t$. 
\end{itemize}
Thus every edge in $W$ has head different from $t$,
and we claim that $S'\cup T$ and $S'\cup \overleftarrow{T}$ are acyclic. For suppose that $C$
is a directed cycle of one of $S'\cup T$ and $S'\cup \overleftarrow{T}$. If $t\notin V(C)$, then $C$ is a directed cycle
of one of $R'\cup (T\setminus t)$ and $R'\cup \overleftarrow{T}\setminus t$, which is impossible. Thus
$t\in V(C)$, and since no edge has head $t$ in $S'\cup \overleftarrow{T}$, it follows that 
$C$ is a directed cycle of $S'\cup T$. We may assume that $C$ is optimal.
Let $e$ be the edge of $C$ with tail $t$, and let its head be $v$. Since $tt_1$
is the only edge of $S'\cup T$ with head $t$, it follows that this edge belongs to $C$; and so there is a directed path
of $(S'\cup T)\setminus t$ from $v$ to $t_1$. Since $R'\cup (T\setminus t)$ is acyclic, there is no edge
in $R'$ with tail $t_1$ and head $v$; consequently there is no edge of $R$ in $W'$ with ends $t_1,v$; in particular,
$e\notin E(R)$; and so $T(v,t_1)$ is a directed path, from the definition of $R$.
Since both ends of this path belong to $V(C)$, (1) implies that $T(v,t_1)\subseteq C$, and so $T(v,t)\subseteq C$; but
$T(v,t)$ is not a directed path, a contradiction. This proves that $S'\cup T$ and $S'\cup \overleftarrow{T}$ are acyclic,
and so the result holds if $d=1$. This proves (2).

\bigskip
Assign an orientation ``clockwise'' to $\Sigma$. We may assume that $W$ is nonempty; for if it is empty, we may replace it
by $\{e\}$ for any edge $e$ of $S$, and if $E(S)=\emptyset$ there is nothing to prove. 
Since $W$ is a $t$-wedge containing no edges of $T$, and $d\ge 2$, 
we can choose two other $t$-wedges $W_1,W_2$, such
that $W,W_1,W_2$ are pairwise disjoint and their union contains all edges of $S\cup T^-$ incident with $t$, and moreover
$W_1,W_2$ both contain at least one edge in $T^-$. We may assume that for every choice of $w\in W, w_1\in W_1$ and $w_2\in W_2$,
the three edges $w,w_1,w_2$ are in clockwise order around $t$.

For $i = 1,2$, let $T_i$ be the subdigraph of $T$ formed by the union of all paths of $T$ with one end $t$ and containing
some edge in $W_i$. Thus, $T_1,T_2$ are directed trees; their union is $T$, and they only have the vertex $t$ in common.
For $i = 1,2$, let $N_i$ be the set of all vertices $v$ of $T_i\setminus \{t\}$ such that
$T(t,v)$ is a directed path; and let $M_i = V(T_i)\setminus (N_i\cup \{t\})$.

For $i = 1,2$, let $S_i$ be
the subgraph of $S$ with vertex set $V(S)$ and edge set the set of all edges $uv$ of $S$ such that either
$u,v\in V(T_i)$ or
one of $u,v\in M_i$.
\\
\\
(3) {\em For $i = 1,2$, there is a directing of $S_i$, forming a digraph $S_i'$, with the following properties:
\begin{itemize}
\item for each edge $e = uv$ of $S_i$ with $u\notin V(T_i)$, the head of $e$ is $v$, and there is 
no directed path of  $S_i'\cup T_i$ or of 
$S_i'\cup \overleftarrow{T_i}$  from $v$ to $t$
\item each edge in $W\cap E(S_i)$ has tail $t$
\item the subdigraphs of $S_i'\cup T_i$
and $S_i'\cup \overleftarrow{T_i}$ induced on $V(T_i)$ are acyclic.
\end{itemize}
}
\noindent From the symmetry between $W_1$ and $W_2$ (reversing the orientation of $\Sigma$ if necessary) 
we may assume that $i = 1$ for definiteness. 
Let $Q$ be the graph with vertex set $V(T_1)$ and edge set $E(S_1)$, with incidence relation as follows:
\begin{itemize}
\item for every edge $e=uv$ of $S$ with $u,v\in V(T_1)$, $e$ has ends $u,v$ in $Q$
\item for every edge $e=uv$ of $S$ with $v\in M_1$ and $u\notin V(T_1)$, $e$ has ends $t,v$ in $Q$.
\end{itemize}
Note first that $Q\cup T_1^-$ is planar. A planar drawing can be obtained from the drawing of $S_1\cup T^-$
by contracting the edges of $T_2$. For the same reason, there 
is a $t$-wedge $W_1$ of this
drawing which contains all edges of $W\cap E(S_1)$ and all edges $e$ of $S_1$ such that $e$ is incident in $S$ with a vertex 
not in $V(T_1)$.

Now for every edge $e$ of $Q$, the path $P$ of $T_1$ between its ends is not a directed path of $T_1$ (because let $u,v$
be the ends of $e$ in $S$; if $u,v\in V(T_1)$, then $P=T(u,v)$, and if
$u\notin V(T_1)$ say, then $v\in M_1$, and $P=T(t,v)$ which is not directed).

Since $|V(T_1)|<|V(T)|$, the inductive hypothesis implies that there is a directing of $Q$, forming a digraph $Q'$ say, such that
$Q'\cup T_1, Q'\cup \overleftarrow{T_1}$ are acyclic and every edge in $W_1$ has tail $t$.
The same directing is the desired directing of $S_1$. This proves (3).

\bigskip
A {\em line} means a subset of $\Sigma$ homeomorphic to the closed
interval $[0,1]$, and we define its {\em ends} in the natural way.
We say an edge of $S$ is {\em low} if 
it has one end in $V(T_1)\setminus \{t\}$ and the other in $V(T_2)\setminus \{t\}$ and with at least 
one end in $N_1\cup N_2$. 
\\
\\
(4) {\em Either there is no low edge in $S$, or there is a line $F$ in $\Sigma$ with the following properties:
\begin{itemize}
\item $F$ has ends in $V(T_1)\setminus \{t\}$ and $V(T_2)\setminus \{t\}$ (say $v_1,v_2$ respectively), 
and either $v_1\in N_1$ or $v_2\in N_2$ 
\item the interior of $F$ is disjoint from the drawing of $S\cup T^-$
\item the closed curve formed by the union of $F$ and the path $T(v_1,v_2)$ bounds two closed discs $D_1,D_2$
in $\Sigma$, where all edges in $W$ are drawn in the interior of $D_2$
\item every low edge of $S$ is drawn within $D_1$.
\end{itemize}
}
\noindent Let us take a cyclic ordering of the edges of $S\cup T^-$ incident with $t$, first those in $W$, and then those in $W_1$,
and finally those in $W_2$. Between any two edges that are consecutive in this ordering, there is a region of the drawing,
in the natural sense.
Let $r$ be the region that comes between the last edge of $W_1$ and the first edge of $W_2$. Let $w$ be an edge in $W$.
It follows that for every edge $e$ of $S$ with ends $v_1 \in V(T_1)\setminus \{t\}$ and $v_2 \in V(T_2)\setminus \{t\}$, 
the cycle formed by the path $T(v_1,v_2)$ and the edge $e$ bounds two closed discs in $\Sigma$, one containing $r$ and 
not $w$, and the other containing $w$ and not $r$. We denote the first of these by $D(e)$.

We may assume that there is a low edge in $S$;
let $e$ be such an edge, chosen with $D(e)$ maximal.
Let $e$ have ends $v_1\in V(T_1)\setminus \{t\}$
and $v_2\in V(T_2)\setminus \{t\}$ say. Thus either $v_1\in N_1$ or $v_2\in N_2$. 

We claim that every low edge $f$ is drawn within $D(e)$. Let $f$ have ends 
$u_i\in V(T_i)$ for $i = 1,2$. If $f=e$ our claim is true, so we assume that $f\ne e$.
The intersection of the 
paths $T(u_1,u_2), T(v_1,v_2)$ is a path, and so the boundaries of $D(e),D(f)$ intersect in a path. Thus either
one of $D(e),D(f)$ includes the other, or they intersect in precisely the intersection of their boundaries, 
or they have union $\Sigma$. The last two alternatives are impossible, because $D(e), D(f)$ both contain $r$, and both
are disjoint from $w$. Moreover, $D(f)$ does not properly include $D(e)$, from the choice of $e$; so $D(e)$ includes $D(f)$,
and in particular $f$ is drawn in $D(e)$.

Let $D'(e)$ be the second disc in $\Sigma$ with the same boundary as $D(e)$.
Choose a line $F$ with ends $v_1,v_2$ and with interior in $D'(e)$, such that
the disc in $D'(e)$ bounded by the union of $e$ and $F$ meets the drawing only in $e$. Let
$D_2$ be the disc in $D'(e)$ bounded by $F$ and the path $T(v_1,v_2)$, and let $D_1$ be the second disc in $\Sigma$ bounded
by the same curve. It follows that every low edge is drawn in the interior of $D_1$. This proves (4).
\\
\\
(5) {\em We may assume that there is a line $F$ in $\Sigma$ with the following properties:
\begin{itemize}
\item $F$ has ends in $V(T_1)$ and $N_2$ (say $f_1,f_2$ respectively)
\item the interior of $F$ is disjoint from the drawing of $S\cup T^-$
\item the closed curve formed by the union of $F$ and the path $T(f_1,f_2)$ bounds two closed discs $D_1,D_2$
in $\Sigma$, where all edges in $W$ are drawn in the interior of $D_2$
\item every low edge of $S$ is drawn within $D_1$.
\end{itemize}
}
\noindent If there is a low edge in $S$,
the claim follows from (4), exchanging $W_1,W_2$ and reversing the orientation of $\Sigma$ if necessary.
If there is no low edge, let $f_1 = t$, and let $f_2\in N_2$ be adjacent to $t$ in $T$; then $f_1,f_2$ are on a common region, 
and we may choose $F$ joining $f_1,f_2$
with interior in this region. This proves (5).

\bigskip

\begin{figure} [h!]
\centering
\includegraphics{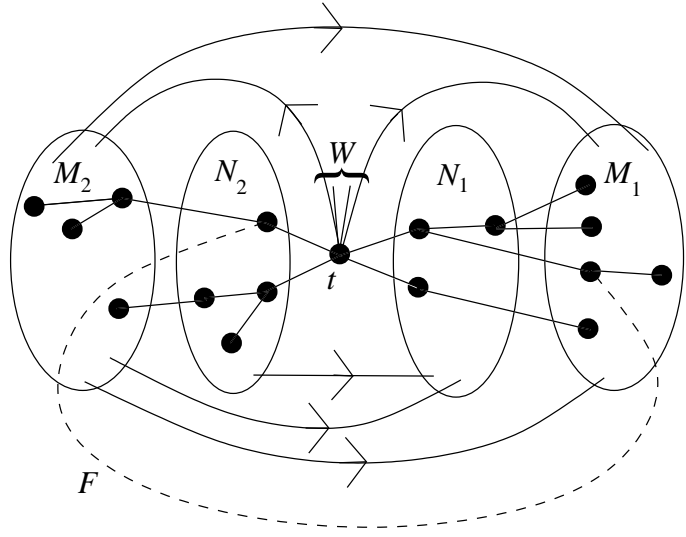}
\caption{Combining solutions. 
$T_1,T_2$ are the subtrees to the right and left of $t$.}
\label{proof}
\end{figure}

Let $f_1,f_2,D_1,D_2$ be as in (5), and let $S_1', S_2'$ be as in (3).
Now we are ready to construct a directing $\eta$ of $S$. (See figure 5.)
For each edge $e$ of $S$:
\begin{itemize}
\item if $e\in E(S_1)$, let $\eta(e)$ be its head in $S_1'$
\item if $e\in E(S_2)\setminus E(S_1)$ and $e$ is drawn in $D_1$, let $\eta(e)$ be the tail of $e$ in $S_2'$
\item if $e\in E(S_2)\setminus E(S_1)$ and $e$ is drawn in $D_2$, let $\eta(e)$ be the head of $e$ in $S_2'$
\item if $e\notin E(S_1)\cup E(S_2)$, then $e$ has ends in $N_1$ and $N_2$; let $\eta(e)$ be its end in $N_1$.
\end{itemize}
We claim that this satisfies the theorem. 
\\
\\
(6) {\em If $e\in W$ then $\eta(e)\ne t$.}
\\
\\
Let $e=tv$. Thus
$v\notin N_1\cup N_2$, because $T(t,v)$ is not a directed path. Consequently if $v\in V(T_1)$ then $e\in E(S_1)$ and so
$\eta(e)\ne t$. If $v\in V(T_2)$, then  $e\in E(S_2)\setminus E(S_1)$; and since $e$ is drawn in $D_2$, 
it follows that $\eta(e)$ is the head of $e$ in $S_2'$,
and therefore not $t$. This proves (6).
\\
\\
(7) {\em $\eta(e) = v_1$ for every edge $e$ of $S$ with ends $v_1\in V(T_1)\setminus \{t\}$ and $v_2\in V(T_2)\setminus \{t\}$.}
\\
\\
For if $v_1\in M_1$ then $e\in E(S_1)$ and hence $\eta(e)$ is its head in $S_1'$; but its head in $S_1'$ is in $V(T_1)$ from
the choice of $S_1'$, and so equals $v_1$. Thus we assume that $v_1\in N_1$, and so $e$ is a low edge, and therefore
drawn within $D_1$. If $e\in E(S_2)$, then $\eta(e)$ is the tail of $e$ in $S_2'$; but $v_2$ is the head of $e$ in $S_2'$, and so
its tail in $S_2'$ is $v_1$ as required. Finally, if $e\notin E(S_1)\cup E(S_2)$, then $\eta(e) = v_1$ from the definition of $\eta$. This
proves (7).

\bigskip
Let $S'$ be the digraph $(S,\eta)$; it remains to show that $S'\cup T$ and $S'\cup \overleftarrow{T}$ are acyclic. 
By reversing all edges of $T$ if necessary, it suffices to prove that $S'\cup T$ is acyclic. Suppose
then that $C$ is a directed cycle of $S'\cup T$, and choose it optimal.
\\
\\
(8) {\em For $j = 1,2$, $V(C)\not\subseteq V(T_j)$.}
\\
\\
For $j = 1$ this is clear, since if $V(C)\subseteq V(T_1)$ then $C$ is a directed cycle of 
$S_1'\cup T_1$, which is impossible. We assume then that $V(C)\subseteq V(T_2)$, and so every edge of $C$
belongs to $S_2\setminus S_1$.
If no edge of $C$ is drawn in the interior of $D_1$, then $C$ is a directed cycle of 
$S_2'\cup T_2$, which is impossible. If no edge of $C$ is drawn in the interior
of $D_2$, then $C$ is a directed cycle of 
$\overleftarrow{S_2'}\cup T_2$; but this is acyclic since
$S_2'\cup \overleftarrow{T_2}$ is acyclic. Thus some edge of $C$ is drawn in the interior of $D_1$,
and some edge in the interior of $D_2$, and so the intersection of $C$ with the boundary of $D_1$ is disconnected. But 
every vertex of $C$ in the boundary of $D_1$ belongs to $T(f_2,t)$, since $V(C)\subseteq V(T_2)$, and this is a directed path
of $T$, contrary to (1). This proves (8).

\bigskip

From (7) and (8), it follows that $t\in V(T)$, and there is a unique edge $e$ of $C$ with one end in 
$V(T_1)\setminus \{t\}$ and one in $V(T_2)\setminus \{t\}$; say the edge $c_1c_2$, where
$c_j\in V(C) \cap V(T_j)\setminus \{t\}$ for $j = 1,2$. By (7) $\eta(e) = c_1$, and since $C$ is a directed cycle
it follows that
$C$ is the union of
\begin{itemize}
\item  a directed path $P_2$ of $S'\cup T$ from $t$ to $c_2$, with $V(P_2)\subseteq V(T_2)$, 
\item and the edge $e$,
and 
\item a directed path $P_1$ of $S'\cup T$ from $c_1$ to $t$, with $V(P_1)\subseteq V(T_1)$. 
\end{itemize}
It follows that $P_1$ is a directed path of $S_1'\cup T_1$, and so $e\notin E(S_1)$ from the choice of $S_1'$. Thus $c_1\in N_1$, and
(1) implies that $P_1 = T(c_1,t)$. Moreover, 
$e$ is a low edge, and consequently is drawn within $D_1$. Since $t\in V(P_2)$, and (1) implies that the intersection
of $P_2$ with $T(f_2,t)$ is connected, it follows that no edge of $C$ is drawn within $D_2$, and in particular
$P_2$ is a directed path of $\overleftarrow{S_2'}\cup T_2$, and hence $\overleftarrow{P_2}$ is a directed path
of $S_2'\cup \overleftarrow{T_2}$ from $c_2$ to $t$. From the choice of $S_2'$, it follows that $e\notin E(S_2)$, and so
$c_2\in N_2$. Since $c_2,t\in V(C)$, (1) implies that $P_2 = T(t,c_2)$; but since $P_1 = T(t,c_1)$, it follows that
$T(c_1,c_2)$ is a directed path (since it is a subpath of $C$), contrary to hypothesis. This proves \ref{wedgethm}.~\bbox

We remark that the proofs of both parts of \ref{orientthm} work by proving a somewhat stronger statement, that a specific
path of $S$ can be made directed (for a $t$-wedge becomes a path under planar duality). Perhaps the following
holds, which would be a common strengthening of \ref{caterpillar}, \ref{wedgethm}, and \ref{biasconj}:

\begin{thm}\label{pathbiasconj}
{\bf Conjecture: }Let $S$ be a tree, and let $\mathcal{B}$ be a bias in $S$, such that $|D_S(X)|\ge 2$ for each $X\in \mathcal{B}$.
Let $P$ be a path of $S$.
Then there is a directing of $S$, forming a digraph $S'$ say, such that
$D^+_{S'}(X),D^-_{S'}(X)\ne\emptyset$ for each $X\in \mathcal{B}$, and such that $P$ becomes a directed path of $S'$.
\end{thm}

\end{document}